\documentclass[11pt,a4paper]{article}

\usepackage{amsmath,amssymb}
\usepackage{verbatim}
\usepackage{enumerate}
\usepackage{hyperref}
\usepackage[T1]{fontenc}
\usepackage[utf8]{inputenc}
\usepackage{authblk}
\usepackage{mathtools}
\usepackage{a4,graphicx}
\usepackage{float}
\usepackage{caption}
\usepackage{subcaption}

\newtheorem{dfn}{Definition}[section]
\newtheorem{lem}[dfn]{Lemma}
\newtheorem{theor}[dfn]{Theorem}

\newtheorem{rem}[dfn]{Remark}
\newtheorem{cor}[dfn]{Corollary}
\newtheorem{ex}[dfn]{Example}

\newcommand{\proof}{\noindent{\bfseries Proof. }}

\newcommand{\thline}{\rule{5.75in}{.01in}}

\title{Fixed points for weak contractions in metric type spaces
}
\date{\vspace{-5ex}}
\author[$\dagger$]{\small Ya\'e Ulrich Gaba }

\affil[$\dagger$]{\small Department of Mathematics and Applied Mathematics, University of Cape Town, Rondebosch 7701, South Africa}

\affil[$\dagger$]{\small gabayae2@gmail.com }

\numberwithin{equation}{section}
\begin{document}

\maketitle

\begin{abstract}
In this article, we prove some fixed point theorems in metric type spaces. This article is just a generalization some results previously proved in \cite{niyi-gaba}. In particular, we give some coupled common fixed points theorems under weak contractions. These results extend well known similar results existing in the literature.
\section*{\small Keyword: Metric Type Spaces; Fixed Point.}
\end{abstract}
\thline

\section{Introduction}
The Banach contraction principle is a fundamental result in fixed point theory. Due to its importance, several authors have obtained many interesting extensions and generalizations, from single fixed point to coupled fixed point. In \cite{niyi-gaba}, Olaniyi et al. gave a brief justification on the origin of metric type space and considered it as a suitable framework to generalize some fixed point theorems. Metric type spaces generalize the notion of metric spaces, even if both structures are very close . In this article, we prove fixed point theorems in metric type spaces. We also give coupled fixed point results. For the convenience of the reader, we shall recall some definitions but for a detailed expos\'e of the definition and examples, the interested reader is referred to \cite{niyi-gaba}.

\section{Prelimaries}

\begin{dfn}
Let $X$ be a nonempty set, and let the function $D:X\times X \to [0,\infty)$ satisfy the following properties:
\begin{itemize}
\item[(D1)] $D(x,x)=0$ for any $x \in X$;
\item[(D2)] $D(x,y)=D(y,x)$ for any $x,y\in X$;
\item[(D3)] $D(x,y) \leq \alpha \big( D(x,z_1)+D(z_1,z_2)+\cdots+D(z_n,y) \big)$ for any points $x,y,z_i\in X,\ i=1,2,\ldots, n$ where $n\geq 1$ is a fixed natural number and some constant $\alpha\geq 1$.
\end{itemize}
The triplet $(X,D,\alpha)$ is called a \textbf{metric type space}.
\end{dfn}

The concepts of Cauchy sequence and convergence for a sequence in a metric type space are defined in the same way as defined for a metric space. For the interested reader the definitions can be obtained in \cite{kam}. Moreover, for $\alpha=1$, we recover the classical  metric, hence metric type generalizes metric. It is worth mentioning that if $(X,D,\alpha)$ is a metric type space, then for any $\beta \geq \alpha$, $(X,D,\beta)$ is also a metric type space. Hence, in the sequel we shall denote $(X,D,\alpha)$ simply as $(X,D)$ when there is no confusion.

\begin{ex}(Compare \cite{niyi-gaba})
Let $X=\{ a,b,c \}$ and the mapping $D:X\times X \to [0,\infty)$ defined by $D(a,b)=1/5,\ D(b,c)=1/4,\ D(a,c)=1/2$, $D(x,x)=0$ for any $x \in X$ and $D(x,y)=D(y,x)$ for any $x,y\in X$.
Since 
\[
\frac{1}{2} = D(a,c)> D(a,b)+D(b,c)=\frac{9}{20},
\]
\end{ex}
then we conclude that $X$ is not a metric space space. Nevertheless, with $\alpha=2$, it is very easy to check that $(X,D,2)$ is a metric type space.

\begin{dfn}
A subset $S$ of a metric type space $(X,D,\alpha)$ is said to be \textbf{totally bounded} if given $\varepsilon>0$ there exists a finite set of points $\{s_1,s_2,\cdots,s_n \} \subseteq S$, called $\varepsilon$-net, such that given any $s\in S$ there exists $i \in \{ 1,2,\cdots,n\}$ for which $D(s,s_i)\leq \varepsilon$.
\end{dfn}

\begin{dfn}
Let $(X,D,\alpha)$ be a metric type space. If for any sequence $(x_n)$ in $X$, there is a subsequence $(x_{n_k})$ of $(x_{n})$ such that $(x_{n_k})$ is convergent in $X$. Then $X$ is called a {\bf sequentially compact} metric type space.
\end{dfn}

\vspace*{0.5cm}

Some topological properties concerning metric type spaces can be read in \cite{niyi-gaba}. We recall some of these results, as we will make use of them.

\begin{lem}(Compare \cite{niyi-gaba})
Every metric type space $(X,D,\alpha)$ is Hausdorrf.
\end{lem}

\begin{lem}(Compare \cite{niyi-gaba})
A Cauchy sequence $(x_n)$ in a metric type space $(X,D,\alpha)$ is always bounded.
\end{lem}

\begin{lem}(Compare \cite{niyi-gaba})
Let $(x_n)$ be a Cauchy sequence in a metric type space $(X,D,\alpha)$. Then $(x_n)$ converges to $x$ if and only if it has a subsequence that converges to $x$.
\end{lem}

\section{First Fixed Point Theorems}
In this section, we prove some fixed point results in metric type space. In particular we generalize the contractive conditions in the literature by replacing the constants with functions. First, we state the following useful lemma.

\begin{lem}(Compare \cite{niyi-gaba})\label{cochi1}
Let $(y_n)$ be a sequence in a metric type space $(X,D,\alpha)$ such that 
\begin{equation}\label{cochi}
D(y_n,y_{n+1}) \leq \lambda D(y_{n-1},y_n)
\end{equation}
for some $\lambda>0$ with $\lambda <\frac{1}{\alpha}$.
Then $(y_n)$ is Cauchy.
\end{lem}

We are now in a position to state the main fixed point theorem of this section.

\begin{theor}\label{theorem1}
 Let $(X,D,\alpha)$ be a complete metric type space and let $f: X \to X$ be a continuous function. Suppose that there exist functions $\eta, \lambda,\zeta,\mu,\xi:X \to [0,1)$ which satisfy the following for $x,y \in X:$ 
\begin{enumerate}
\item[(1)] $g(f(x)) \leq g(x)$ whenever $g \in \{\eta, \lambda,\zeta,\mu,\xi\} ;$ 
\item[(2)]  $\eta(x)+\lambda(x)+\zeta(x)+\mu(x)+2\alpha \xi(x)<1;$
\item[(3)]  $D(f(x),f(y)) \leq H(x) * V(x,y)$ where $H(x)$ and $V(x,y)$ are vectors in $\mathbb{R}^5$  

with $H(x):=[\eta(x),\lambda(x),\zeta(x),\mu(x),\xi(x)]$, 

$V(x,y):= [D(x,y),D(x,f(x)),D(y,f(y)),D(f(x),y),D(x,f(y))]$ and $*$ the usual inner product on $\mathbb{R}^5$.

\end{enumerate}
Then, $f$ has a unique fixed point.
\end{theor}

\proof
Let $x_0 \in X$ be arbitrary and fixed, and we consider the sequence $x_{n+1}=fx_n=f^{n+1}x_0, $ for all $n\in \mathbb{N}$. If we take $x=x_{n-1}$ and $y=x_n$ in $(3)$ we have 
\begin{align*}
 D(x_n,x_{n+1}) & = D(fx_{n-1},fx_n)\\
                & \leq H(x_{n-1}) * V(x_{n-1},x_n) \\
                & = H(fx_{n-2}) * V(x_{n-1},x_n).\\
 \end{align*}
Since for all $n \in \mathbb{N},$  $V(x_{n-1},x_n)=[D(x_{n-1},x_n),D(x_{n-1},x_n)D(x_{n},x_{n+1}),0,D(x_{n-1},x_{n+1})$, the above becomes
\begin{align*}
 D(x_n,x_{n+1}) & \leq H(fx_{n-2}) * V(x_{n-1},x_n)\\
                & \leq H(x_{n-2}) * V(x_{n-1},x_n) \\
                & \leq \eta(x_{n-2})D(x_{n-1},x_n)+\lambda(x_{n-2})D(x_{n-1},x_n)+\zeta(x_{n-2})D(x_n,x_{n+1})\\
                & + \alpha \xi(x_{n-2})[D(x_{n-1},x_n)+D(x_n,x_{n+1})] \\
                & \vdots \\
                & \leq \eta(x_0)D(x_{n-1},x_n)+\lambda(x_0)D(x_{n-1},x_n)+\zeta(x_0)D(x_n,x_{n+1})\\
                & + \alpha \xi(x_0)[D(x_{n-1},x_n)+D(x_n,x_{n+1})]. \\
 \end{align*}
So 
$$D(x_n,x_{n+1}) \leq \frac{\eta(x_0)+\lambda(x_0)+\alpha \xi(x_0)}{1-\zeta(x_0)-\alpha \xi(x_0)}D(x_{n-1},x_n) .  $$

Thus, by Lemma \ref{cochi1}, $(x_n)$ is a Cauchy sequence in $X$. Because of completeness of $X$ and continuity of $f$, there exists $x^* \in X$ such that $x_n \to x^*$ and $x_{n+1}=f(x_n)\to f(x)$. Since $X$ is Hausdorff, $f(x^*)=x^*$ .

\vspace*{0.5cm}

\underline{ Uniqueness} If $y^*$ is a fixed point of $f$, then
\begin{align*}
D(x^*,y^*)= D(fx^*,fy^*) & \leq  \eta(x^*)D(x^*,y^*) + \lambda(x^*)D(x^*,fx^*)  \\
                         & + \zeta(x^*)D(y^*,fy^*) + \mu(x^*) D(fx^*,y^*) \\
                         & + \xi(x^*)D(x^*,fy^*) \\
                         & =   \eta(x^*)D(x^*,y^*) + \mu(x^*) D(x^*,y^*) + \xi(x^*)D(x^*,y^*) \\
                         & = (\eta(x^*) + \mu(x^*) + \xi(x^*))D(x^*,y^*)).
\end{align*}
Therefore $D(x^*,y^*)=0$ i.e. $x^*=y^*$ since $\eta(x^*) + \mu(x^*) + \xi(x^*)<1$.

\begin{cor}\label{cor1}
Let $(X,D,\alpha)$ be a complete metric type space and let $f: X \to X$ be a continuous function. Suppose that there exist functions $\eta, \lambda,\zeta,\mu,\xi:X \to [0,1)$ which satisfy the following for $x,y \in X:$ 
\begin{enumerate}
\item[(1)] $g(f(x)) \leq g(x)$ whenever $g \in \{\eta, \lambda,\zeta,\mu,\xi\} ;$ 
\item[(2)]  $\eta(x)+2\lambda(x)+2(1+\alpha)\mu(x)<1;$
\item[(3)]  $D(f(x),f(y)) \leq \eta(x)D(x,y)+\lambda(x)(D(x,f(x))+D(y,f(y)))+\\
\mu(x)(D(f(x),y)+D(x,f(y))) $. 
\end{enumerate}
Then, $f$ has a unique fixed point.
\end{cor}

\proof
We can prove this result by applying Theorem \ref{theorem1} to $\lambda(x)=\zeta(x)$ and $\mu(x)=\xi(x)$.

\begin{cor}\label{cor2}
Let $(X,D,\alpha)$ be a complete metric type space and let $f: X \to X$ be a continuous function. Suppose that there exist constants $a,\beta,\gamma,k,l\geq 0$ which satisfy the following for $x,y \in X:$ 

$D(f(x),f(y)) \leq aD(x,y)+ \beta D(x,f(x))+ \gamma D(y,f(y))+kD(f(x),y)+lD(x,f(y))$ 
and $a,\beta,\gamma,k,l\geq 0$ with $a+\beta+\gamma+k+2\alpha l< 1$.
Then, $f$ has a unique fixed point.
\end{cor}

\proof
We can prove this result by applying Theorem \ref{theorem1} to $ H(x)=(a,\beta,\gamma,k,l)$.

\begin{cor}\label{cor3}
Let $(X,D,\alpha)$ be a complete metric type space and let $f: X \to X$ be a continuous function and

$$D(f(x),f(y)) \leq  \beta [D(x,f(x))+  D(y,f(y))]$$ 
for all $x,y \in X$  and $\beta\in [0,1/2)$. Then, $f$ has a unique fixed point.
\end{cor}

\begin{cor}\label{cor4}
Let $(X,D,\alpha)$ be a complete metric type space and let $f: X \to X$ be a continuous function and

$$D(f(x),f(y)) \leq  \beta [D(f(x),y)+  D(x,f(y))]$$ 
for all $x,y \in X$  and $\beta\in [0,1/(1+2\alpha))$. Then, $f$ has a unique fixed point.
\end{cor}

\begin{cor}\label{cor5}
Let $(X,D,\alpha)$ be a complete metric type space and let $f: X \to X$ be a continuous function and

$$D(f(x),f(y)) \leq  aD(x,y)+ \beta D(f(x),y)$$ 
for all $x,y \in X$  and $a+ \beta < 1$. Then, $f$ has a unique fixed point.
\end{cor}

The next corollary is a generalization of Banach contraction principle.
\begin{cor}\label{cor6}
Let $(X,D,\alpha)$ be a complete metric type space and let $f: X \to X$ be a continuous function and 

$$D(f(x),f(y)) \leq  aD(x,y)$$ 
for all $x,y \in X$  and $a < 1$. 
\end{cor}

We conclude this section by he following theorem.

\begin{theor}\label{theorem2}
Let $(X,D,\alpha)$ be a complete metric type space and let $f: X \to X$ be a continuous function and for all $x,y \in X$

\begin{equation}\label{star}
\varepsilon(x,y)  \leq s D(x,y)  +t D(x,f^2(x))   \tag{*}            
\end{equation}
with $$\varepsilon(x,y)=a D(f(x),f(y)) + \beta D(x,f(x))+ \gamma D(y,f(y))  
 + k D(x,f(y)) + lD(y,f(x)),$$
and where $s\geq a\geq l,$ $\gamma\geq k\geq t,$ $a+k>0$ and $0\leq (s-l)/(a+k)<1/\alpha$. 

Then, $f$ has a fixed point. Moreover, if $\alpha>1$, the fixed point is unique.
\end{theor}

\proof
For simplicity, we set $fx:=f(x)$ and subsequently. 
Let $x_0 \in X$ be arbitrary and fixed, and we consider the sequence $x_{n+1}=fx_n=f^{n+1}x_0, $ for all $n\in \mathbb{N}$. If we take $x=x_{n-1}$ and $y=x_n$ in \eqref{star} we have 

\begin{align}
a D(f(x_{n-1}),f(x_{n})) &+ \beta D(x_{n-1},f(x_{n-1}))+ \gamma D(x_{n},f(x_{n})) \nonumber      \\ 
 &+ k D(x_{n-1},f(x_{n})) + lD(x_{n},f(x_{n-1}))  \nonumber    \\ 
      &\leq s D(x_{n-1},x_{n})  +t D(x_{n-1},f^2(x_{n-1}))     
\end{align}
Rewriting his inequality as

\begin{align}
a D(x_{n},x_{n+1}) &+ \beta D(x_{n-1},x_{n})+ \gamma D(x_{n},x_{n+1}) \nonumber      \\ 
 &+ k D(x_{n-1},x_{n+1}) + l D(x_{n},x_{n})   \nonumber   \\ 
      &\leq s D(x_{n-1},x_{n})  +t D(x_{n-1},x_{n+1})     
\end{align}
implies that

\begin{equation}\label{alter}
(a+\gamma)D(x_{n},x_{n+1}) +(k-t)D(x_{n-1},x_{n+1}) \leq (s-\beta) D(x_{n-1},x_{n}) .
\end{equation}
Since $k\geq t$, we have 
\[
(a+\gamma)D(x_{n},x_{n+1})  \leq (s-\beta) D(x_{n-1},x_{n}), 
\]
which entails that 
\[
D(x_{n},x_{n+1})  \leq  \frac{s-\beta}{a+\gamma}  D(x_{n-1},x_{n}) \leq \frac{s-l}{a+k}  D(x_{n-1},x_{n}), 
\]

Thus, by Lemma \ref{cochi1}, $(x_n)$ is a Cauchy sequence in $X$. Because of completeness of $X$ and continuity of $f$, there exists $x^* \in X$ such that $x_n \to x^*$ and $x_{n+1}=f(x_n)\to f(x)$. Since $X$ is Hausdorff, $f(x^*)=x^*$ .

\vspace*{0.5cm}

\underline{Uniqueness} We assume here that $K>1$. If $y^*$ is a fixed point of $f$, then for $x=x*$  and $y=y^*$ in \eqref{star}, we obtain

\begin{align}
a D(f(x^*),f(y^*)) &+ \beta D(x^*,f(x^*))+ \gamma D(y^*,f(y^*))  
 + k D(x^*,f(y^*)) + lD(y^*,f(x^*))  \tag{*} \\ 
      &\leq s D(x^*,y^*)  +t D(x^*,f^2(x^*)).      \nonumber         
\end{align}
Hence,

\[
(a+k)D(x^*,y^*) +l D(x^*,y^*) \leq s D(x^*,y^*) \Longleftrightarrow D(x^*,y^*) \leq \frac{s-l}{a+k} D(x^*,y^*)
\]

which is absurd except $D(x^*,y^*)=0$, i.e. $x^*=y^*$.

\begin{rem}
Observe in \eqref{alter} that the term $D(x_{n-1},x_{n+1})$ can be split via the triangle inequality. This yields to 
\[
(a+\gamma)D(x_{n},x_{n+1}) \leq (k-t) \alpha [D(x_{n-1},x_{n})+ D(x_n, x_{n+1}) ] + (s-\beta) D(x_{n-1},x_{n}),
\]
i.e.
\[
[a+\gamma + (k-t)\alpha ]  D(x_n, x_{n+1}) \leq [(k-t) \alpha + s-\beta] D(x_{n-1},x_{n})
\]
or 
\[
D(x_n, x_{n+1}) \leq \frac{[(k-t) \alpha + s-\beta]}{[a+\gamma -(k-t)\alpha ] } D(x_{n-1},x_{n}).
\]
In this case, we shall require that $\frac{[(k-t) \alpha + s-\beta]}{[a+\gamma -(k-t)\alpha ] } < 1/\alpha$ which implies that $(s-l)/(a+k)<1/\alpha$.
Hence there is no loss of generality if we do not use the triangle inequality at this stage.
\end{rem}

The next sections present a generalization of the results proved above, namely, we remplace the set of one variable functions $\eta, \lambda,\zeta,\mu,\xi:X \to [0,1)$ by a set of functions of two variables $\eta, \lambda,\zeta,\mu,\xi:X\times X \to [0,1)$ to which we add an additional condition.

\section{Fixed point of generalized contraction mappings}

\begin{theor}\label{res1}
Let $(X,D,K)$ be a complete metric type space and let $T: X \to X$ a self mapping on $X$ such that for each $x,y \in X$

\begin{align}\label{Equation1}
D(Tx,Ty)\leq \ &\alpha(x,y) D(x,y) + \beta(x,y)D(x,Tx)+ \gamma(x,y)D(y,Ty) \nonumber \\
&+\delta(x,y)[D(x,Ty)+D(y,Tx)]
\end{align}
where $\alpha,\beta,\gamma,\delta$ are functions from $X\times X \to [0,1)$ such that 
\begin{equation}\label{Equation2}
\lambda = \sup \{\alpha(x,y)+\beta(x,y)+\gamma(x,y) +  2K\delta(x,y): x,y  \in X\}<1,
\end{equation}
then 
$T$ has a unique fixed point, say $x^*$, the orbit $\{T^nx\}$ converges to the fixed point $x^*$ for any $x\in X$. Moreover $x^*$ is such that for each $x\in X$
\[
D(T^nx,x^*)\ \leq \frac{\lambda^n}{1-\lambda}D(x,Tx).
\]
\end{theor}
\proof
Let $x \in X$. Let $(x_n)$ be the sequence defined by $x_0=x, \ x_1=Tx_0, \ \cdots x_{n+1}=Tx_n$. From \eqref{Equation1},
\begin{align*}
D(Tx_n,Tx_{n+1})= D(Tx_{n-1},Tx_n)\leq \ &\alpha D(x_{n-1},x_{n}) + \beta D(x_{n-1},x_{n}) + \gamma D(x_{n},x_{n+1}) \nonumber \\
&+\delta [D(x_{n-1},x_{n+1})+D(x_{n},x_n)],
\end{align*}
or 
\begin{align}\label{Equation3}
D(Tx_n,Tx_{n+1})= D(Tx_{n-1},Tx_n)\leq \ &\alpha D(x_{n-1},x_{n}) + \beta D(x_{n-1},x_{n}) + \gamma D(x_{n},x_{n+1}) \nonumber \\
&+\delta D(x_{n-1},x_{n+1}).
\end{align}
Hence, using the triangle inequality, we have

\begin{align}\label{Equation4}
D(x_{n-1},x_{n+1})& \leq K [D(x_{n-1},x_{n+1})+D(x_{n-1},x_{n+1}) ] \nonumber \\
                  &\leq 2K \max \{ D(x_{n-1},x_{n}),D(x_{n},x_{n+1})\}.
\end{align}
By \eqref{Equation4}, equation \eqref{Equation3} turns to be

\begin{align*}
D(x_{n},x_{n+1})& \leq (\alpha+ \beta+\gamma) \max \{ D(x_{n-1},x_{n+1}),D(x_{n-1},x_{n+1})\}\\
                  &+ 2K \max \{ D(x_{n-1},x_{n}),D(x_{n},x_{n+1})\}.
\end{align*}
Then
\[
D(x_{n},x_{n+1}) \leq \lambda \max \{ D(x_{n-1},x_{n}),D(x_{n},x_{n+1})\}.
\]
Since $\lambda<1,$ then
\begin{equation}
D(x_{n},x_{n+1}) \leq \lambda  D(x_{n-1},x_{n}).
\end{equation}

By inductivity, we obtain that 
\[
D(x_{n},x_{n+1}) \leq \lambda  D(x_{n-1},x_{n}) \leq \lambda.\lambda D(x_{n-2},x_{n-1})\leq \cdots \leq \lambda^n D(x,Tx).
\]
By triangle inequality, for $m>n$ we get

\begin{align*}
D(x_n,x_m)&\leq K[D(x_{n},x_{n+1})+D(x_{n+1},x_{n+2})+\cdots + D(x_{m-1},x_{m})]\\
          &\leq K[\lambda^n +\lambda^{n+1}+\cdots +\lambda^{m-1} ]D (x,Tx)\\
          & \frac{\lambda^n}{1-\lambda}KD (x,Tx)
\end{align*}
or
\begin{equation}\label{Equation7}
D(x_n,x_m)\frac{\lambda^n}{1-\lambda}KD (x,Tx)
\end{equation}
Letting $m,n \to\infty$, in \eqref{Equation7} implies that $(x_n)$ is a Cauchy sequence. Since $(X,D,K)$ is complete, there exists $x^*\in X$ such that 
\begin{equation}\label{Equation8}
\underset{n\to \infty}{\lim}x_n=x^*.
\end{equation}
We now show that $T(x^*)=x^*$. Form the contractive condition, we know that
\begin{align*}
D(Tx^*,Tx_n) & \leq  \alpha D(x^*,x_n) + \beta D(x^*,Tx^*) +\gamma D(x_n,Tx_n)\\
              & + \delta [D(u,Tx_n)+D(x_n,Tu)] \\
              & \leq \lambda \max \{  D(x^*,x_n) , D(x^*,Tx^*),D(x_n,x_{n+1}),
              D(u,x_{n+1}),D(x_n,Tu) \}.
\end{align*}

As $n\to \infty$, we obtain
\begin{equation}
D(Tx^*,x^*)\leq \lambda D(x^*,Tx^*).
\end{equation}
Since $\lambda<1$, we conclude that $T(x^*)=x^*$ because $DT(x^*,x^*)=0$.

For uniqueness, assume $y^* \in $ is a fixed point of $T$. Then, we have

\begin{align*}
D(x^*,y^*)&=D(Tx^*,Ty^*)\\
          &  \leq \alpha D(x^*,y^*) + \beta D(x^*,Tx^*)+ \gamma D(y^*,Ty^*) \nonumber \\
&+\delta[D(x^*,Ty^*)+D(y^*,Tx^*)] \\
           & \leq (\alpha +2\delta)D(x^*,y^*) \leq \lambda D(x^*,y^*),
\end{align*}
and $D(x^*,y^*)=0$ since $\lambda<1$.

Taking the limit in \eqref{Equation7} as $m\to \infty$, we get $D(T^nx,x^*)\ \leq \frac{\lambda^n}{1-\lambda}D(x,Tx)$ for each $n$. The proof is complete.

Mappings which satisfy \eqref{Equation1} and \eqref{Equation2} are called {\bf generalized contraction}. From the proof of Theorem \ref{res1}, we have 

\begin{cor}
Let $(X,D,K)$ be a complete metric type space and let $T: X \to X$ self mapping on $X$ such that for each $x,y \in X$
\begin{equation}\label{Equation10}
D(Tx,Ty)\leq \lambda \max\left\lbrace D(x,y),D(x,Tx),D(y,Ty),\frac{1}{2}[D(x,Ty)+D(y,Tx)]\right\rbrace
\end{equation}
where $\lambda \in (0,1)$. Then $T$ has a unique fixed point.
\end{cor}

\begin{cor}
Let $(X,D,K)$ be a complete metric type space and let $T: X \to X$ a self mapping on $X$ such that for each $x,y \in X$

\begin{equation}\label{Equation1m}
D(Tx,Ty)\leq  \lambda_1 D(x,Tx)+ \lambda_2 D(y,Ty)+ \lambda_3 D(x,Ty)+\lambda_4 D(y,Tx)
\end{equation}
where $\lambda_i \in [0,1)$ and $\lambda_1 +\lambda_2+K(\lambda_3+\lambda_4)< \min\{ 1,2/K\}, i=1,2,3,4$. Then $T$ has a unique fixed point.
\end{cor}

\section{Common fixed point of generalized contraction mappings}

Given a non empty set $X$ and $\{T_\alpha\}_{\alpha\in J}$ a family of selfmappings on $X$ where $J$ is an index set. A point $u\in X$ is called common fixed point for the family $\{T_\alpha\}_{\alpha\in J}$ if and only if $u=T_\alpha(u)$ for each $\alpha$.

\begin{theor}
Let $(X,D,K)$ be a complete metric type space and $\{T_\alpha\}_{\alpha\in J}$ a family of selfmappings on $X$. If there exists a fixed $\beta\in J$ such that for each $x,y \in X$ such that for each $\alpha\in J$:

\begin{equation}\label{Equation18}
D(T_\alpha x,T_\beta y)\leq  \lambda \max \left\lbrace D(x,y),D(x,T_\alpha x),D(y,T_\beta y),\frac{1}{2}[ D(x,T_\beta y), D(y,T_\alpha x)]\right\rbrace
\end{equation}
for some $\lambda=\lambda(\alpha) \in [0,1)$ with $\lambda K<1$ and for each $x,y \in X$. Then all the $T_\alpha$ have a unique common fixed point, which is unique for each $T_\alpha, \alpha\in J$.
\end{theor}

\proof
Let $\alpha \in J$ and $x\in X$ be arbitrary. We consider the sequence defined by $x_0=x,\ x_{2n+1}=T_\alpha x_{2n}, \ x_{2n+2}=T_\beta x_{2n+1}$. Form \eqref{Equation18} we get 
\begin{align*}
D(x_{2n+1}, x_{2n+2}) &= D(T_\alpha x_{2n},T_\beta x_{2n+1})\\
                      & \leq \lambda \max
                  \left\lbrace
                \begin{array}{ll}
                  D(x_{2n}, x_{2n+1}), D(x_{2n}, x_{2n+1}), D(x_{2n+1}, x_{2n+2}),\\
                  \frac{1}{2} [D(x_{2n}, x_{2n+2})+D(x_{2n+1}, x_{2n+1})]
                \end{array}
              \right\rbrace
\end{align*}

Since 
\begin{equation}\label{Equation19}
D(x_{2n}, x_{2n+2}) \leq K [ D(x_{2n}, x_{2n+1})+D(x_{2n+1}, x_{2n+2})]
\end{equation}
so, 
\begin{align*}
 \frac{1}{2} D(x_{2n}, x_{2n+2}) & \leq K [D(x_{2n}, x_{2n+1})+D(x_{2n+1}, x_{2n+2}) ] \\
                                 & \leq K \max\{ D(x_{2n}, x_{2n+1}),D(x_{2n+1}, x_{2n+2})\},
\end{align*}
we have 
\[
D(x_{2n+1}, x_{2n+2}) \leq K \max\{ D(x_{2n}, x_{2n+1}),D(x_{2n+1}, x_{2n+2})\}.
\]
Hence, as $\lambda K<1$

\[
D(x_{2n+1}, x_{2n+2}) \leq \lambda K D(x_{2n}, x_{2n+1}).
\]

Similarly, we get that $D(x_{2n}, x_{2n+1}) \leq \lambda K D(x_{2n-1}, x_{2n}).$ thus for any $n\geq 1$ we have 

\begin{equation}\label{Equation20}
D(x_{n},x_{n+1}) \leq \lambda K D(x_{n-1},x_{n})\leq K^2 D(x_{n-2},x_{n-1})\leq \cdots \leq (\lambda K)^{n} D(x_{0},x_{1})  
\end{equation}

From \eqref{Equation20} and the triangle inequality, for $m>n$ we get
\begin{align*}
D(x_n,x_m) & \leq K[D(x_{n},x_{n+1})+D(x_{n+1},x_{n+2})+\cdots + D(x_{m-1},x_{m})]\\
           & \leq K (  (\lambda K)^{n} D(x_{0},x_{1})  +(\lambda K)^{n+1} D(x_{0},x_{1})  + \cdots + (\lambda K)^{m-1} D(x_{0},x_{1}) ) \\
           & \leq K \left[  \lambda K)^{n}  +(\lambda K)^{n+1}  + \cdots + (\lambda K)^{m-1}\right] D(x_{0},x_{1}) \\
          & \leq K \frac{(\lambda K)^n}{1-(\lambda K)} D(x_{0},x_{1}). 
\end{align*}
If we let $n \to \infty$, we conclude that $(x_n)$ is a Cauchy sequence in $X$. Because of completeness of $X$ there exists $x^* \in X$ such that $x_n \to x^*$.

From \eqref{Equation18}, we have

\begin{align*}
D(T_\beta x^*,x_{2n+1}) &= D(T_\beta x^*,T_\beta x_{2n})\\
                      & \leq \lambda \max
                  \left\lbrace
                \begin{array}{ll}
                  D(x^*, x_{2n}), D(x^*, T_\beta x^*), D(x_{2n}, x_{2n+1}),\\
                  \frac{1}{2} [D(x^*, x_{2n+1})+D(x_{2n}, T_\beta x^*)]
                \end{array}
              \right\rbrace .
\end{align*}
Taking the limit as $n\to \infty$, we get $D(T_\beta x^*, x^*)\leq \lambda D(T_\beta x^*, x^*)$. Therefore $D(T_\beta x^*, x^*)=0$, i.e. $T_\beta x^*= x^*$.

Moreover, from \eqref{Equation18} with $x=y=T_\beta x^*= x^*$, we have
\[
D(T_\alpha x^*, x^*) = D(T_\alpha x^*,T_\beta x^*) \leq \lambda(\alpha) \max \left\lbrace D(T_\alpha x^*, x^*), \frac{1}{2}D(T_\alpha x^*, x^*)\right\rbrace
\]
and hence $T_\alpha x^*= x^*$ for all $\alpha \in J$. Thus, all $T_\alpha$ have a common fixed point. Suppose that $w$ is a fixed point of $T_\beta$, it follows as above , that $w$ is a common fixed point for all $T_\alpha$. Thus, from \eqref{Equation18} we have $D(z,x^*)=D(T_\beta x^*, T_\alpha w)\leq \lambda D(x^*,w)$ and so $x^*=w$. Thus, $x^*$ is the unique common fixed point for all $\{T_\alpha\}_{\alpha\in J}$. This completes the proof.

\section{Weak contractions}
In this section, we present some common fixed point under contractive conditions for $w$-compatible mappings in metric type spaces.

We start with the following definition.

\begin{dfn}(Compare \cite{ce})
Let $X$ be a non empty set. An element $(x,y) \in X\times X$ is called:

\begin{enumerate}
\item[(E1)] a \textbf{coupled fixed point} of the mapping $F:X\times X \to X$ if $F(x,y)=x$ and $F(y,x)=y$;

\item[(E2)] a \textbf{coupled coincidence point} of the mappings $F:X\times X \to X$ and $g:X \to X$ if $F(x,y)=g(x)$ and $F(y,x)=g(y)$. In this case $(gx,gy)$ is called the \textbf{coupled point of coincidence};

\item[(E3)] a \textbf{coupled common fixed point} of the mappings $F:X\times X \to X$ and $g:X \to X$ if $F(x,y)=gx=x$ and $F(y,x)=gy=y$.
\end{enumerate}
\end{dfn}

\begin{dfn}
Let $X$ be a non empty set. The mappings $F:X\times X \to X$ and $g:X \to X$ are said to be \textbf{$w$-compatible} if $g(F(x,y))=F(gx,gy)$ whenever $F(x,y)=g(x)$ and $F(y,x)=g(y)$. 

Observe here that if $(x,y)$ is a coupled common fixed point of the mappings $F:X\times X \to X$ and $g:X \to X$, then likewise $(y,x)$.
\end{dfn}

\begin{theor}\label{thm2.2}
Let $(X,D,K)$ a metric type space and $F:X\times X \to X$ and $g:X \to X$ be maps such that for all $x,y,u,v \in X$

\begin{align}\label{Eq2.1}
D(F(x,y),F(u,v)) & \leq \alpha_1 D(gx,gu) + \alpha_2 D(F(x,y),gx) +\alpha_3 D(gy,gv) \nonumber \\
                 & + \alpha_4 D(F(u,v),gu) +\alpha_5 D(F(x,y),gu) \nonumber \\ 
                 & + \alpha_6 D(F(u,v),gx),
\end{align}
where $\alpha_i$ for $i=1,\cdots ,6$ are nonnegative constants with 
\begin{equation}\label{Eq2.2}
2K(\alpha_1+\alpha_3)+(K+1)(\alpha_2+\alpha_4)+(K^2+K)(\alpha_5+\alpha_6)<2.
\end{equation}
If $F(X\times X)\subset g(X) $ and $g(X)$ is complete, then $F$ and $g$ have a coupled coincidence point in $X$.
\end{theor}

\proof
Let $x_0,y_0 \in X$ be arbitrary and set 
\[
g(x_1)=F(x_0,y_0),\ g(y_1)=F(y_0,x_0),\cdots , \ g(x_{n+1})= F(x_n,y_n),\ g(y_{n+1})=F(y_n,x_n)
\]
This construction is always possible since $F(X\times X)\subset g(X) $. From \eqref{Eq2.1}, we have
\begin{align}\label{Eq2.3}
D(gx_{n+1},gx_n) & = D(F(x_n,y_n),F(x_{n-1},y_{n-1})) \nonumber \\
                 & \leq \alpha_1 D(gx_n,gx_{n-1}) + \alpha_2 D(F(x_n,y_n),gx_n) +\alpha_3 D(gy_n,gy_{n-1}) \nonumber \\
                 & \ \ \ + \alpha_4 D(F(x_{n-1},y_{n-1}),gx_{n-1}) +\alpha_5 D(F(x_n,y_n),gx_{n-1}) \nonumber \\
                 &\ \ \ + \alpha_6 D(F(x_{n-1},y_{n-1}),gx_n) \nonumber \\
                 &  \leq \alpha_1 D(gx_n,gx_{n-1}) + \alpha_2 D(g(x_{n+1}),gx_n) +\alpha_3 D(gy_n,gy_{n-1}) \nonumber \\
                 & \ \ \ + D(gx_n,gx_{n-1}) +\alpha_5 D(gx_{n+1},gx_{n-1})+ \alpha_6 D(gx_n,gx_n) \nonumber \\
                 &  \leq \alpha_1 D(gx_n,gx_{n-1}) + \alpha_2 D(g(x_{n+1}),gx_n) +\alpha_3 D(gy_n,gy_{n-1}) \nonumber \\
                 & \ \ \ + D(gx_n,gx_{n-1}) + K \alpha_5 [D(gx_{n+1},gx_n)+ D(gx_n,gx_{n-1})].         
\end{align} 
It follows that
\begin{equation}\label{Eq2.4}
(1-\alpha_2-K\alpha_5)D(gx_{n+1},gx_n) \leq (\alpha_1+\alpha_4+K\alpha_5)D(gx_n,gx_{n-1}) + \alpha_3 D(gy_n,gy_{n-1}). 
\end{equation}

Similarly, we have 
\begin{equation}\label{Eq2.5}
(1-\alpha_2-K\alpha_5)D(gy_{n+1},gy_n) \leq (\alpha_1+\alpha_4+K\alpha_5)D(gy_n,gy_{n-1}) + \alpha_3 D(gx_n,gx_{n-1}). 
\end{equation}

From 
\[
D(gx_{n+1},gx_n) = D(F(x_n,y_n),F(x_{n-1},y_{n-1}))= D(F(x_{n-1},y_{n-1}),F(x_n,y_n)) = D(gx_n,gx_{n+1}), 
\]
and applying \eqref{Eq2.1} to $D(gx_n,gx_{n+1})$, we derive that 

\begin{equation}\label{Eq2.7}
(1-\alpha_4-K\alpha_6)D(gx_n,gx_{n+1}) \leq (\alpha_2+\alpha_4+K\alpha_5)D(gx_{n-1},gx_{n}) +\alpha_3 D(gy_{n-1},gy_n). 
\end{equation}
 
and 
\begin{equation}\label{Eq2.8}
(1-\alpha_4-K\alpha_6)D(gy_n,gy_{n+1}) \leq (\alpha_2+\alpha_4+K\alpha_5)D(gy_{n-1},gy_{n}) +\alpha_3 D(gx_{n-1},gx_n). 
\end{equation}

Now, by adding up \eqref{Eq2.4} and \eqref{Eq2.5} and setting $D_n = D(gx_n,gx_{n+1}) + D(gy_n,gy_{n+1})$, we obtain
\begin{equation}\label{Eq2.9}
(1-\alpha_2-K\alpha_5) D_n \leq (\alpha_1+\alpha_3+\alpha_4+K\alpha_5)D_{n-1}. 
\end{equation}

Similarly, adding \eqref{Eq2.1} and \eqref{Eq2.8} and setting $D_n = D(gx_n,gx_{n+1}) + D(gy_n,gy_{n+1})$, we obtain
\begin{equation}\label{Eq2.10}
(1-\alpha_4-K\alpha_6) D_n \leq (\alpha_1+\alpha_2+\alpha_3+K\alpha_6)D_{n-1}. 
\end{equation}
Finally, adding up \eqref{Eq2.9} and \eqref{Eq2.10}, we have
\begin{equation}\label{Eq2.11}
[2-\alpha_2-\alpha_4-K(\alpha_5+\alpha_6)] D_n \leq [2\alpha_1+\alpha_2+2\alpha_3+\alpha_4+ K(\alpha_5+\alpha_6)]D_{n-1}.
\end{equation}

Hence, for all $n,$
\begin{equation}\label{Eq2.12}
0\leq D_n \leq \lambda D_{n-1} \leq \lambda^2 D_{n-2}\leq \cdots \leq \lambda^n D_0,
\end{equation}
where 

\begin{equation}\label{Eq2.13}
\lambda = \frac{2\alpha_1+\alpha_2+2\alpha_3+\alpha_4+ K(\alpha_5+\alpha_6)}{2-\alpha_2-\alpha_4-K(\alpha_5+\alpha_6)}< \frac{1}{K}.
\end{equation}

If $D_0=0$ then $(x_0,y_0)$ is a coupled coincidence of $F$ and $g$.

Now, assume that $D_0>0$. From Lemma \ref{cochi1}, and since $\lambda< 1/K$, we conclude that $(gx_n)$ and $(gy_n)$ are Cauchy sequences in $X$. Since $g(X)$ is complete, there exixt $x^*,y^* \in X$ such that $gx_n\to gx^*$ and $gy_n\to gy^*$.

Now, we prove that $F(x^*,y^*)=gx^*$ and $F(y^*,x^*)=gy^*$. Indeed, using \eqref{Eq2.1}

\begin{align}\label{Eq2.16}
D(F(x^*,y^*),gx^*) & \leq K [ D(F(x^*,y^*),gx_{n+1})+ D(gx_{n+1},gx^*)  ] \nonumber \\
                   & = [ D(F(x^*,y^*),F(x_n,y_n))+ D(gx_{n+1},gx^*)  ] \nonumber \\
                   & \leq K [  \alpha_1 D(gx,gx_n)+\alpha_2 D(F(x^*,y^*),gx^*)+\alpha_3 D(gy^*,gy_n) \nonumber \\
                   & \ \ + K \alpha_4[D(gx_{n+1},gx^*)+D(gx^*,gx_n)] \nonumber \\
                   & \ \ + K \alpha_5 [ D(x^*,y^*)+D(gx^*,gx_n)]  \nonumber \\
                   & \ \  + \alpha_6 D(gx_{n+1},gx^*) +   D(gx_{n+1},gx^*)                              ].
\end{align}
Therefore
\begin{align}\label{Eq2.17}
D(F(x^*,y^*),gx^*) & \leq \frac{K\alpha_1+K^2(\alpha_4+\alpha_5)}{1-K\alpha_2-K^2\alpha_5}  \nonumber \\
                   & \ \ + \frac{K+K^2\alpha_4 + K \alpha_6}{1-K\alpha_2-K^2\alpha_5}  \nonumber \\
                   & \ \ + \frac{K\alpha_3}{1-K\alpha_2-K^2\alpha_5}.  
\end{align}

Since $gx_n\to gx^*$ and $gy_n\to gy^*$, we have $D(F(x^*,y^*),gx^*)=0$, i.e. $F(x^*,y^*)=gx^*$. Similarly, we can get $F(y^*,x^*)=gy^*$. Therefore $(x^*,y^*)$ is a coupled coincidence point for $F$ and $g$.

\begin{theor}\label{thm2.3}
Let $(X,D,K)$ a metric type space and $F:X\times X \to X$ and $g:X \to X$ be maps which satisfy the conditions of Theorem \ref{thm2.2}. If $F$ and $g$ are $w$-compatible, then $F$ and $g$ have a unique coupled common fixed point, which belongs to the diagonal of $X$.
\end{theor}

\proof

Step 1:

In a first time, we establish the uniqueness of the coupled point of coincidence. Suppose by the way of contraction that $F(x^*,y^*)=gx^*$, $F(y^*,x^*)=gy^*$, $F(a^*,b^*)=ga^*$, $F(b^*,a^*)=gb^*$ for some $x^*,y^*,a^*,b^* \in X$. From \eqref{Eq2.1}, we have 
\[
D(ga^*,gx^*)= D(F(a^*,b^*),F(x^*,y^*) ) \leq (\alpha_1+\alpha_5+\alpha_6)D(ga^*,gx^*)+\alpha_3 D(gb^*,gy^*)
\]
and 
\[
D(gb^*,gy^*)= D(F(b^*,a^*),F(y^*,x^*) ) \leq (\alpha_1+\alpha_5+\alpha_6)D(gb^*,gy^*)+\alpha_3 
D(ga^*,gx^*).
\]
Adding up the above two inequalities, we get 
\[
D(ga^*,gx^*)+ D(gb^*,gy^*) \leq (\alpha_1+\alpha_3+\alpha_5+\alpha_6)[D(ga^*,gx^*)+ D(gb^*,gy^*)]
\]

Since $2K(\alpha_1+\alpha_3)+(K+1)(\alpha_2+\alpha_4)+(K^2+K)(\alpha_5+\alpha_6)<2$, we have $D(ga^*,gx^*)+ D(gb^*,gy^*)=0$, i.e. $ ga^*=gx^*$ and $gb^*=gy^*$. It is also very easy to observe that a similar computation as the one did above establishes that $ga^*=gy^*$ and $gb^*=gx^*$. Thus $ga^*=gb^*$ and $(gx^*,gx^*)$ is the unique coupled point of coincidence of $F$ and $g$. 

Step 2:

Now let $g(x^*)=z$. Then we have $z=g(x^*)=F(x^*,x^*)$. By $w$-compatibility of $F$ and $g$, we have 
\[
g(z)=g(g(x^*))=g(F(x^*,x^*))=F(g(x^*),g(x^*))=F(z,z).
\]

Thus $(gz,gz)$ is a coupled point of coincidence of $F$ and $g$. Therefore $z=gz=F(z,z)$. Consequently $(z,z)$ is the unique coupled common fixed point of $F$ and $g$.

\begin{cor}\label{cor2.4}
Let $(X,D,K)$ a metric type space and $F:X\times X \to X$ and $g:X \to X$ be maps such that for all $x,y,u,v \in X$

\begin{align}\label{Eq2.21}
D(F(x,y),F(u,v)) & \leq \alpha [ D(gx,gu) +  D(F(x,y),gx) ] +  \nonumber \\
                 & \beta [ D(gy,gv) D(F(u,v),gu) ]  \nonumber \\ 
                 & + \gamma [ D(F(x,y),gu) D(F(u,v),gx)],
\end{align}
where $\alpha, \beta$ and $\gamma$ are nonnegative constants with 
\begin{equation}\label{Eq2.22}
(3K+1+)(\alpha+\beta)+2(K^2+K)\gamma <2.
\end{equation}
If $F(X\times X)\subset g(X) $ and $g(X)$ is complete, then $F$ and $g$ have a coupled coincidence point in $X$. Also if $F$ and $g$ are $w$-compatible, then $F$ and $g$ have a unique coupled common fixed point, which belongs to the diagonal of $X$.
\end{cor}

\proof
It follows from Theorems \ref{thm2.2} and \ref{thm2.3} by setting $\alpha_1=\alpha_2=\alpha,$  $\alpha_3=\alpha_4=\beta$ and $\alpha_5=\alpha_6=\gamma$.

\begin{cor}\label{cor2.5}
Let $(X,D,K)$ a metric type space and $F:X\times X \to X$ and $g:X \to X$ be maps such that for all $x,y,u,v \in X$

\begin{equation}\label{Eq2.23}
D(F(x,y),F(u,v))  \leq \alpha D(gx,gu) +  \beta D(gy,gv)  \nonumber \\
\end{equation}
where $\alpha, \beta$  are nonnegative constants with $\alpha+\beta <1/K.$

If $F(X\times X)\subset g(X) $ and $g(X)$ is complete, then $F$ and $g$ have a coupled coincidence point in $X$. Also if $F$ and $g$ are $w$-compatible, then $F$ and $g$ have a unique coupled common fixed point, which belongs to the diagonal of $X$.
\end{cor}

\proof
It follows from Theorems \ref{thm2.2} and \ref{thm2.3} by setting $\alpha_1=\alpha,$  $\alpha_3=\beta$ and $\alpha_2=\alpha_4=\alpha_5=\alpha_6=0$.

\begin{cor}\label{cor2.6}
Let $(X,D,K)$ a metric type space and $F:X\times X \to X$ and $g:X \to X$ be maps such that for all $x,y,u,v \in X$

\begin{equation}\label{Eq2.24}
D(F(x,y),F(u,v)) \leq \alpha D(F(x,y),gu) + \beta  D(F(u,v),gx)
\end{equation}
where $\alpha, \beta$ are nonnegative constants with $\alpha+\beta <2/(K^2+K).$
$(3K+1+)(\alpha+\beta)+2(K^2+K)\gamma <2.$
If $F(X\times X)\subset g(X) $ and $g(X)$ is complete, then $F$ and $g$ have a coupled coincidence point in $X$. Also if $F$ and $g$ are $w$-compatible, then $F$ and $g$ have a unique coupled common fixed point, which belongs to the diagonal of $X$.
\end{cor}

\proof
It follows from Theorems \ref{thm2.2} and \ref{thm2.3} by setting $\alpha_5=\alpha,$  $\alpha_6=\beta$ and $\alpha_1=\alpha_2=\alpha_3=\alpha_4=0$.

We give here an example to illustrate the above results.

\begin{ex}
Let $X=[0,1]$ and define on $X\times X$ the map $D$ given by $D(x,y)=|x-y|^2$. The reader can convince himself that $D$ does not satisfy the triangle inequality and that $(X,D,2)$ is a metric type space. For that, we could use the Minkowski inequality $|x-z|^2 \leq 2 (|x-y|^2+ |y-z|^2)$. Define the map $F:X\times X \to X$ by $F(x,y)=(x+y)/4$ and $g(x)=x$. For $\alpha=\beta=1/8$, it is easy to check that $F$ and $g$ satisfy the condition \eqref{Eq2.23} and that $\alpha+\beta = 1/4 \in [0,1/K)=[0,1/2)$, i.e. 
\begin{equation*}
D(F(x,y),F(u,v)) \leq \frac{1}{8} [D(x,u) + D(y,v)].
\end{equation*}
We apply Corollary \ref{cor2.5} to deduce the existence of coupled fixed point for $F$ and $g$, which is the present case is $(0,0)$.
\end{ex}

\end{document}